\documentclass[12pt]{amsart}
\usepackage{graphicx}
\usepackage{amssymb}
\usepackage{amsfonts}
\usepackage{amsmath,amscd}
\usepackage{array}
\usepackage{rotating}
\usepackage{epsfig}
\usepackage{txfonts}

\newtheorem{example}{Example}[section]

\newtheorem{theorem}[example]{Theorem}

\newtheorem{corollary}[example]{Corollary}

\newtheorem{proposition}[example]{Proposition}

\def\Proof{\noindent \it Proof -- \rm}
\def\qed{\hspace{3.5mm} \hfill \vbox{\hrule height 3pt depth 2 pt width 2mm}
\bigskip}

\def\Co{{\rm Co}}
\def\SP{{\rm SP}}
\def\SW{{\rm SW}}
\def\ol{\bar}
\def\QSym{{\it QSym}}          
\def\Sym{{\bf Sym}}            
\def\NCSF{{\bf Sym}}           
\def\FQSym{{\bf FQSym}}        
\def\FSym{{\bf FSym}}          
\def\WQSym{{\bf WQSym}}        

\def\WSym{{\bf WSym}}          
\def\X{{\bf X}}

\def\Mw{{\bf M}}               

\def\p{{\mathfrak p}}

\def\std{{\rm std}}     

\def\<{\langle}
\def\>{\rangle}
\def\Z{{\bf Z}}
\def\NN{{\mathbb N}}    
\def\KK{{\mathbb K}\, } 

\def\F{{\bf F}}         
\def\SG{{\mathfrak S}}  

\def\shuff#1#2{\mathbin{
\hbox{\vbox{ \hbox{\vrule \hskip#2 \vrule height#1 width 0pt
}%
\hrule}%
\vbox{ \hbox{\vrule \hskip#2 \vrule height#1 width 0pt
\vrule }%
\hrule}%
}}}

\def\shuf{{\mathchoice{\shuff{7pt}{3.5pt}}%
{\shuff{6pt}{3pt}}%
{\shuff{4pt}{2pt}}%
{\shuff{3pt}{1.5pt}}}}%
\def\shuffle{\,\shuf\,}
\def\ad{{\rm ad}}

\title[Binary shuffle bases for quasi-symmetric functions]%
{Binary shuffle bases for quasi-symmetric functions}

\author[J.-C.~Novelli and J.-Y.~Thibon]%
{Jean-Christophe Novelli and Jean-Yves Thibon}

\address[] {Universit\'e Paris-Est Marne-la-Vall\'ee, Laboratoire
d'Informatique Gaspard Monge, 5 Boulevard Descartes \\Champs-sur-Marne \\77454
Marne-la-Vall\'ee cedex 2 \\ FRANCE}
\email[Jean-Christophe Novelli]{novelli@univ-mlv.fr}
\email[Jean-Yves Thibon]{jyt@univ-mlv.fr} 
\date{}

\begin{document}

\begin{abstract}
We construct bases of quasi-symmetric functions whose product rule
is given by the shuffle of binary words, as for multiple zeta values in
their integral representations, and then extend the construction to the
algebra of free quasi-symmetric functions colored by positive integers.
As a consequence, we show that the fractions introduced
in [Guo and Xie, Ramanujan Jour. 25 (2011) 307-317] provide a realization
of this algebra by rational moulds extending that of free quasi-symmetric
functions given in 
[Chapoton et al., Int. Math. Res. Not. IMRN 2008, no. 9, Art. ID rnn018].
\end{abstract}

\maketitle


\section{Introduction}

The algebra of Quasi-symmetric functions $\QSym(X)$ \cite{Ge} is the linear
span of the expressions
\begin{equation}
M_I(X) :=
   \sum_{n_1<n_2<\ldots<n_r} x_{n_1}^{i_1}x_{n_2}^{i_2}\cdots x_{n_r}^{i_r},
\end{equation}
called quasi-monomial functions. Here, $X=\{x_i\}$ is a totally ordered set
of mutually commuting variables and $I=(i_1,\ldots,i_r)$ is a finite sequence
of positive integers (a composition of degree $n=i_1+\cdots+i_r$ of $M_I$).

The $M_I$ are partially symmetric functions (actually, the invariants of
a special action of the symmetric group \cite{Hiv}), and the point is that
$\QSym(X)$ is actually an algebra. The product rule involves an operation on
compositions often called quasi-shuffle \cite{Hof}, or stuffle \cite{BBBL}, or
augmented shuffle \cite{Hiv}.
Let $\sqcupplus$ be defined as follows \cite{TU}: 
\begin{equation}
au \sqcupplus bv = a(u\sqcupplus bv)+b(au\sqcupplus v)+
                 (a\dotplus b)(u\sqcupplus v)\
\end{equation}
\begin{equation} 
a \sqcupplus \epsilon= \epsilon\sqcupplus a = a\,,
\end{equation}
where $a, b\in\NN^*$, $\epsilon$ is the empty word,
$a\dotplus b$ denotes the sum of integers,
 and $u$, $v$ are two compositions regarded as 
words on the alphabet $\NN^*$.

\noindent
Then, given two compositions $I$ and $J$, we have:
\begin{equation}\label{quasish}
M_I M_J = \sum_K \< K\ |\ I\sqcupplus J \> M_K
\end{equation}
where the notation $\<x|y\>$ means ``the coefficient of
$x$ in $y$''.
For example:
\begin{equation}
\begin{split}
M_{2 1}M_{1 2}=&\,
2 M_{{2112}}+M_{{2121}}+M_{{213}}+M_{{222}}\\
&+M_{{1212}}+2 M_{{1221}}+M_{{123}}+M_{{141}}+M_{{312}}+M_{{321}}+M_{{33}}\,.
\end{split}
\end{equation}

As is well known, this quasi-shuffle, which contains all the terms
of the ordinary shuffle $I\shuffle J$ together with contractions
obtained by adding two consecutive parts coming from different terms, defines
a structure actually isomorphic to an ordinary shuffle algebra. The best way
to see this is to recall that $\QSym$ is a Hopf algebra, and that its dual is
the algebra $\Sym$ of noncommutative symmetric functions \cite{MR,NCSF1}. 
Since $\Sym$
can be freely generated by a sequence containing one primitive element $P_n$  in
each degree, the dual $Q_I$ of the multiplicative basis
$P^I=P_{i_1}\cdots P_{i_r}$ satisfies
\begin{equation}\label{shufcomp}
Q_I Q_J = \sum_K \< K\ |\ I\shuffle J \> Q_K\,.
\end{equation}
Thus, the quasi-shuffle algebra over the positive integers is isomorphic to
the shuffle algebra over the same set, and it is straightforward to see that
this is true in general \cite{Hof}.

This isomorphism plays an important role in the theory of multiple zeta values
(MZV), also called polyzetas or Euler-Zagier sums \cite{Zag}. Indeed, these
are specializations of the quasi-monomial functions obtained by setting
\begin{equation}
x_n = \frac1n\quad\text{or}\quad x_{-n}=\frac1n\quad \text{(there are two
conventions)}.
\end{equation}
The second convention being  more frequent, we shall set
\begin{equation}
\zeta(I)=\sum_{n_1>n_2>\ldots>n_r}\frac1{n_1^{i_1}n_2^{i_2}\cdots n_r^{i_r}}
\end{equation}
so that the convergent ones are those with $i_1>1$.

Hence, the convergent MZV satisfy the product rule \eqref{quasish}.
However, one of the big mysteries of the theory relies on the fact
that they also satisfy another product formula, the shuffle relation,
which is {\it not} the shuffle of compositions as in \eqref{shufcomp}.
This formula comes from the following integral representation.
Let
\begin{equation}
\omega_0(t)=\frac{dt}{t}\quad\text{and}\quad \omega_1(t)=\frac{dt}{1-t}\,.
\end{equation}
Then, for a composition $I$ of $n$,
\begin{equation}\label{eqintZ}
\zeta(I)
= \int_{0<t_1<t_2\ldots<t_n}
     \omega_{\epsilon_1}(t_n)\wedge\omega_{\epsilon_2}(t_{n-1})
     \wedge\cdots\wedge\omega_{\epsilon_n}(t_1)
\end{equation}
where $\epsilon$ is the word $\epsilon_I=0^{i_1-1}10^{i_2-1}1\cdots
0^{i_r-1}1$.
Then, a standard property of iterated integrals (often referred to as Chen's
lemma) implies that
\begin{equation}
\zeta(I)\zeta(J)=\sum_K \<\epsilon_K|\epsilon_I\shuffle \epsilon_J\>\zeta(K)\,.
\end{equation}

The existence of such a formula for a specialization of quasi-symmetric
functions raises the question of the existence of a similar one for the
generic case.
That is, does there exist a basis $Z_I$ of $\QSym$ such that

\begin{equation}\label{shufZ}
Z_I Z_J=\sum_K \<\epsilon_K|\epsilon_I\shuffle \epsilon_J\>Z_K\,?
\end{equation}

We shall answer this question in three different ways. First, a counting
argument shows easily that such bases do exist. This does not however provide
a pratical way to construct them. Next, we propose a recusive algorithm
allowing to construct such a basis from any basis $Q_I$ satisfying
\eqref{shufcomp}.  Finally, we describe another construction and we obtain an
explicit combinatorial formula for matrices expressing the dual basis $X_I$ of
$Z_I$ over the dual basis $P_I$ of any basis $Q_I$ as above. 

The interesting point of this last construction is that it involves in a
crucial way the Hopf algebra of set partitions $\WSym$ (symmetric functions in
noncommuting variables), and its relations with the commutative combinatorial
Hopf algebras exhibited in \cite{HNT-co}. Indeed, our matrices are obtained
from two simple statistics of set partitions, which are in some way compatible
with the Hopf structure. The row sums, for example, yield a known form of the
noncommutative Bell polynomials, and the column sums correspond to new analogs
of these. 

Finally, we extend our construction to $\FQSym^{(\NN^*)}$, the algebra of
colored free quasi-symmetric functions, with the set of positive integers as
color set \cite{NT-col}. Namely, we exhibit bases whose product rule coincides
with that of the MZV fractions introduced by Guo and Xin \cite{GX}, thus
providing a realization of this algebra by rational moulds extending that of
$\FQSym$ given in \cite{CHNT}.

\section{Lyndon words and Lyndon compositions}

We shall work in the shuffle algebra $\KK_{\shuffle}\<a,b\>$ over two letters
$a$ and $b$ (rather than 0,1), with $a<b$, $\KK$ being a field of
characteristic 0.
We then replace the notation $\epsilon_I$ by $W_I$, so that
\emph{e.g.}, $W_{213}=abbaab$.

By Radford's theorem \cite{Rad}, $\KK_{\shuffle}\<a,b\>$ is a polynomial
algebra freely generated by Lyndon words. Apart from $a$, all
Lyndon words over $a,b$ end by a $b$.
Words ending by $b$ are in bijection with compositions, and Lyndon words
ending by $b$ correspond to Lyndon compositions. Precisely, $I$ is anti-Lyndon
(a Lyndon composition for the opposite order on the integers)
iff $W_I$ (or $\epsilon_I$) is a Lyndon word.

For example, the Lyndon words of length 6 are
\begin{equation}
aaaaab,\ aaaabb,\ aaabab,\ aaabbb,\ aababb,\ aabbab,\ aabbbb,\ ababbb\,
abbbbb,
\end{equation}
and their encodings by compositions are
\begin{equation}
6,\ 51,\ 42,\ 411,\ 321,\ 312,\ 3111,\ 2211,\ 21111.
\end{equation}

Since $\QSym$ is the shuffle algebra over the positive integers, it is
a polynomial algebra over (anti-) Lyndon compositions. Hence it is
isomorphic to the subalgebra $K_{\shuffle}\<a,b\>b$ (spanned by
words ending by $b$) of the shuffle algebra over two letters.
We then have

\begin{proposition}
There exists a basis $Z_I$ of $\QSym$ satisfying Eq. \eqref{shufZ}:
\begin{equation*}
Z_I Z_J=\sum_K \<\epsilon_K|\epsilon_I\shuffle \epsilon_J\>Z_K \,.
\end{equation*}
\end{proposition}

This argument does not yet give a systematic procedure to build 
a basis with the required properties. An algorithm will be
described in the forthcoming section.

\section{An algorithmic construction}\label{sec:algo}

Our problem is clearly equivalent to the following one:
build a basis $Y_I$ of $K_{\shuffle}\<a,b\>b$ such that

\begin{equation}\label{shufY}
Y_I\shuffle  Y_J = \sum_K \< K\ |\ I\shuffle J \> Y_K\,.
\end{equation}

For an anti-Lyndon composition $L$, we set $Y_L=W_L$ (a Lyndon
word). Thus, we start with
\begin{equation}
Y_1=b,\quad Y_2=ab,\quad
Y_3 = aab,\quad Y_{21} = abb,  \quad Y_{211}=abbb,\ldots
\end{equation}
and applying iteratively \eqref{shufY}, we obtain
\begin{align}
Y_1\shuffle Y_1&= 2Y_{11} = b\shuffle b =2bb,\\
Y_1\shuffle Y_2& = Y_{12}+Y_{21} = b\shuffle ab = bab + 2abb,\\
Y_1\shuffle Y_{21}&= Y_{121}+2Y_{211},
\end{align} 
and so on, from which we deduce
\begin{equation}
Y_{11}=W_{11}, \ Y_{12}=W_{12}+W_{21}, \ Y_{121}=W_{121}+W_{211},\ldots 
\end{equation}
This leads to a triangular system of equations, which determines each
$Y_I$ as a linear combination with nonnegative integer coefficients
of the $W_J$, such that $\ell(J)=\ell(I)$ and $J\ge I$ for the lexicographic
order on compositions.

Here are the first transition matrices (entry $(I,J)$ is the coefficient of
$W_I$ in $Y_J)$. The indexation is the same all over the paper: the
compositions are sorted according to the reverse lexicographic order.
Note that zeroes have been replaced by dots to enhance readibility.

With $n=4$, the indexations of the non-trivial blocks are respectively
$(31),(22),(13)$ and $(211), (121), (112)$ and the matrices are

\begin{equation}
\label{n4k23}
\left(\begin{array}{rrr}
1 & 2 & 1 \\
. & 1 & 1 \\
. & . & 1
\end{array}\right)
\qquad
\qquad
\left(\begin{array}{rrr}
1 & 1 & 1 \\
. & 1 & 1 \\
. & . & 1
\end{array}\right)
\end{equation}

With $n=5$, the indexations of the non-trivial blocks are respectively
$(41),(32),(23),(14)$ and $(311), (221), (212), (131), (122), (113)$
and the matrices are

\begin{equation}
\label{n5k23}
\left(\begin{array}{rrrr}
1 & . & 6 & 1 \\
. & 1 & 2 & 1 \\
. & . & 1 & 1 \\
. & . & . & 1
\end{array}\right)
\qquad
\qquad
\left(\begin{array}{rrrrrr}
1 & . & 6 & 1 & . & 1 \\
. & 1 & 1 & 1 & 2 & 1 \\
. & . & 1 & . & 1 & 1 \\
. & . & . & 1 & 2 & 1 \\
. & . & . & . & 1 & 1 \\
. & . & . & . & . & 1
\end{array}\right)
\end{equation}

One may observe the following properties of these matrices:

\begin{itemize}
\item They are triangular and block diagonal if compositions
are ordered by reverse length-lexicographic order 
(e.g., $4, 31,22,13, 211,121,112,1111$ for $n=4$).
\item The block for length $k$ has therefore dimension $\binom{n-1}{k-1}$.

\item Moreover, the sum of the entries of this block is the
Stirling number of the second kind $S(n,k)$ (the number of
set partitions of an $n$-set into $k$ blocks.)
\end{itemize}

The first two properties are obvious. 
To prove the third one, define
\begin{equation}
P_{n,k} =\sum_{I\vDash n,\ \ell(I)=k}Y_I.
\end{equation}
This can be rewritten as a sum over partitions
\begin{equation}
P_{n,k}=\sum_{\mu=(1^{m_1}2^{m_2}\cdots p^{m_p})\vdash n,\ \ell(\mu)=k}
Y_{1^{m_1}}\shuffle Y_{2^{m_2}}\shuffle\cdots\shuffle Y_{p^{m_p}},
\end{equation}
and since
\begin{equation}
Y_{i^{m_i}} = \frac{W_i^{\shuffle m_i}}{m_i!}
\end{equation}
we see that $P_{n,k}$ is the coefficient of $x^ny^k$
in
\begin{equation}
\exp_{\shuffle}\left\{ \sum_{m\ge 1} x^my W_m\right\}.
\end{equation}
Applying the character $W_n\mapsto 1/n!$ of $\KK_{\shuffle}\<a,b\>b$,
we recognize the generating series of Stirling numbers of the second kind.

We shall now construct ``better'' matrices, sharing all these properties, and
for which a closed formula can be given.

\section{Dual approach}

\subsection{Stalactic conguences and word-symmetric functions}

Alternatively, we can try to build a basis $X_I$ of $\Sym$ whose
coproduct is the binary unshuffle dual to the product of the $Z_I$.
As above, we start with a basis $V_I$ of $\Sym$ whose coproduct
is dual to the shuffle of compositions. For example, we can
take $V_I=X^I$, where $X_n$ is a sequence of primitive generators.
But there are other choices. For example, in \cite{HNT-co}, such
a basis is obtained from a realization of $\Sym$ as a quotient
of $\WSym$, a Hopf algebra based on set partitions.

By definition, $\WSym(A)$  is the subspace of $\KK\<A\>$ spanned by the
orbits of 
the symmetric group
$\SG(A)$ acting on $A^*$ by automorphisms.
These orbits are naturally labelled by set partitions of $[n]$, the orbit
corresponding to a partition $\pi$ being constituted of the words
\begin{equation}
w = a_1\ldots a_n,
\end{equation}
such that $a_i=a_j$ iff $i$ and $j$ are in the same block of $\pi$.
The sum of these words will be denoted by $\Mw_\pi$.

For example,
\begin{equation}
\Mw_{\{\{1,3,6\},\{2\},\{4,5\}\}} := \sum_{a\not=b ; b\not=c; a\not=c}
abacca.
\end{equation}

It is known that the natural coproduct of $\WSym$ (given as usual by the
disjoint union of mutually commuting alphabets) 
is cocommutative~\cite{BRRZ} and that $\WSym$ is
free over connected set partitions.

There are several ways to read an integer composition from a set partition.
First, one can order the blocks according to the values of their minimal
or maximal elements, and record the lenghts of the blocks. For example,
we can order
\begin{equation}
\pi= \left\{\left\{3, 4\right\}, \left\{5\right\}, \left\{1, 2, 6\right\}\right\}
\end{equation} 
in two ways, obtaining two set compositions
\begin{equation}
\Pi' = (126|34|5)\quad \text{and}\quad \Pi''=(34|5|126)
\end{equation}
and the integer compositions
\begin{equation}
K'(\pi)= (321)\quad \text{and}\quad K''(\pi)=(213)\,.
\end{equation}

In \cite{HNT-co}, it is proved that  the two-sided
ideal of $\WSym$ generated by the differences $\Mw_\pi-\Mw_{\pi'}$ for $K'(\pi)=K'(\pi')$
is a Hopf ideal and that the quotient is  isomorphic to $\Sym$.
This quotient can be interpreted in terms of a congruence on the free
monoid $A^*$, the (left) {\it stalactic congruence}, introduced in
\cite{HNT-co}. It is generated by the relations
\begin{equation}
awa \equiv' aaw\quad \text{for $a\in A$ and $w\in A^*$.}
\end{equation}
Thus, each word $w$ is congruent to the word $P(w)$ obtained by moving each 
letter towards its leftmost occurence. Recording the original
position of each letter by a set partition $Q(w)$, we obtain a bijection
which is formally similar to (although much simpler than) the Robinson-Schensted
correspondence, and $\WSym$ can be characterized as the costalactic algebra
is the same way as $\FSym$ is the coplactic algebra~\cite{NCSF6}.
One can of course also define a right stalactic congruence
by
\begin{equation}
awa \equiv'' waa\quad \text{for $a\in A$ and $w\in A^*$.}
\end{equation}
%
The same properties hold with the symmetric condition $K''(\pi)=K''(\pi')$,
which amounts to take the right stalactic quotient.
From now on, we write $K$ instead of $K''$, and denote by $V_I$ the equivalence
class of $\Mw_\pi$ such that $K(\pi)=I$.

One can also record the values of the minimal and maximal elements of the blocks
in the form of a composition. The minimal elements form a subset of $[n]$
always containing $1$, so that we can decrement each of them and remove 0 so as to obtain
a subset of $[n-1]$ which can be encoded by a composition $C'(\pi)$ of $n$.
Similarly, the maximal elements form a subset always containing $n$, so that removing
$n$, we obtain again a subset of $[n-1]$ and a composition $C''(\pi)$ of $n$.
For example, with $\pi$ as above,
\begin{equation}
C'(\pi)= (222)\quad \text{and}\quad C''(\pi)=(411)\,.
\end{equation}
We shall choose the second option and write $C(\pi)$ for $C''(\pi)$.

Clearly, both $C(\pi)$ and $K(\pi)$ have the same length, which is the
number of blocks of $\pi$.

\subsection{A sub-coalgebra of $\WSym$}

Aside from being a quotient of $\WSym$, $\NCSF$ is also a sub coalgebra of
$\WQSym$. More precisely,

\begin{theorem}
The sums
\begin{equation}
\X_J = \sum_{C(\pi)=J}\Mw_\pi
\end{equation}
span a sub-coalgebra of $\WSym$, and
\begin{equation}
\Delta\X_J = \sum_{K,L}\<W_J|W_K\shuffle W_L\>\X_K\otimes \X_L\,.
\end{equation}
\end{theorem}

\Proof
This can be seen on an appropriate encoding.

Note first that the partitions of a set $S$ of integers can be encoded by the
set $\SW(S)$ of signed words whose letters are the elements of $S$, signed
values appearing in increasing order of their absolute values,  and such that
$|w_i|<|w_{i+1}|$ if $w_i$ is unsigned.
Indeed, order the parts according to their maximal elements, sort each part in
increasing order,  read all parts successively and sign (overline) the last
element of each part.
For example,
\begin{equation}
\Pi=(346|5|127) \mapsto 34\ol6\ol512\ol7.
\end{equation}

Now, given a set $S$ and a subset $S'$, define $\SP(S,S')$ as the set of 
set partitions of $S$ whose set of maximal elements of the parts is $S'$.
For example, with $S=\{1,3,5,6\}$ and $S'=\{3,5,6\}$, we have
\begin{equation}\label{eq:w2lsp}
\SP(S,S')= \{ (13|5|6), (3|15|6), (3|5|16)\}.
\end{equation}
Note that $S$ and $S'$ are unambiguously determined by the nondecreasing
word 
$w(S,S')$
whose letters are the elements of $S$, with the elements of $S'$ overlined.
When applied to the signed permutation $\sigma_I$, the identity with the
descents of $I$ overlined, this process encodes the definition of $\X_I$ as a
sum of $\Mw_\pi$.

Finally, given a set partition $\pi$, denote by $\Co(\pi)$ the set of pairs
obtained by splitting the parts of $\pi$ into two subsets in all possible ways.
This encodes the coproduct of an $\Mw_\pi$ expressed in terms of tensor
products  $\Mw_{\pi'}\otimes\Mw_{\pi''}$.

To compute $\Delta\X_I$, we apply successively $\SP$ and $\Co$ to the
signed word $\sigma_I$. 
We then obtain the set $S(I)$ of
pairs of non-intersecting set partitions whose union is a set partition of
$[1,n]$ and whose maximal elements of the blocks are exactly the descents of
$I$.

We shall now check that we get the same result for the coproduct
of a basis $X_I$ of $\Sym$
dual to a basis $Z_I$ of $\QSym$ whose product is given by \eqref{shufZ}.
The coproduct of $X_I$  can be computed by
the following algorithm. Start with the signed permutation $\sigma_I$ defined
as above. It is clearly an encoding of the word $W_I$: an overlined letter
corresponds to a $b$, and the other ones to an $a$. Thus, the coproduct of
such a basis is encoded by the set of pairs of words $(w_1,w_2)$ such that
$\sigma_I$ occurs in $w_1\shuffle w_2$, and such that $w_1$ and $w_2$ end by a
signed letter.
Now, to get this expression from $\Delta\X_I$ 
expanded as a linear combination of terms 
$\Mw_{\pi'}\otimes \Mw_{\pi''}$, one
only needs to apply $\SP$ to each word $w$ separately.
This yields exactly the same set $S(I)$.
Indeed, both sets are multiplicity-free; the second set is a subset of
$S(I)$ since one can only obtain set partitions with given maximal values when
applying $\SP$; and the second set must contain $S(I)$, since any element of
$S(I)$ gives back an unshuffling of $\sigma_I$ when reordering its two subset
partitions.
\qed

\begin{corollary}
For two compositions $I$ and $J$ of $n$, of the same length $k$, let
\begin{equation}
c_{IJ} =
 \#\{\pi\in\Pi_n\ |\ K(\pi)=I\ \text{and}\ C(\pi)=J\}\,,\quad\text{and} \ X_J
 = \sum_I c_{IJ}V_I\,.
\end{equation}
Then, the dual basis $Z_I$ of $X_I$ satisfies the binary shuffle product rule
\eqref{shufZ}.
\end{corollary}

\Proof Denoting by $\bar\Mw_\pi$ the stalactic class of $\Mw_\pi$, we have by definition
\begin{equation}
X_J = \sum_{C(\pi)=J}\bar\Mw_\pi\,.
\end{equation}
Since the canonical projection $p:\ \Mw_\pi\mapsto \bar\Mw_\pi$ is a Hopf
algebra morphism, we have
\begin{equation}
\Delta X_J = \sum_{K,L}\<W_J|W_K\shuffle W_L\> X_K\otimes  X_L\,.
\end{equation}
\qed

\begin{example}{\rm
Let us compute $\Delta(\X_{221})$.
By definition,
\begin{equation}
\X_{211} = \Mw_{12|3|4} + \Mw_{2|13|4} + \Mw_{2|3|14},
\end{equation}
which can be encoded as
\begin{equation}
1\bar2\bar3\bar4 + \bar21\bar3\bar4 + \bar2\bar31\bar4.
\end{equation}
Extracting subwords ending by a signed letter
yields the following list
$S(211)$ of pairs of words:
\begin{equation}
\label{211}
\begin{split}
& (1\bar2\bar3, \bar4), (1\bar2\bar4,\bar3), (1\bar2, \bar3\bar4) \\
& (\bar21\bar3, \bar4), (1\bar3\bar4,\bar2), (1\bar3, \bar2\bar4) \\
& (\bar21\bar4, \bar3), (\bar31\bar4,\bar2), (1\bar4, \bar2\bar3) \\
\end{split}
\end{equation}
together with the symmetrical pairs, where $1$ then belongs to the right part.
The  theorem states that $\Delta \X_{211}$ can be computed
by the unshuffling of the word $1\bar2\bar3\bar4$, which gives the set of words
\begin{equation}
\begin{split}
& (1\bar2\bar3, \bar4), (1\bar2\bar4,\bar3), (1\bar3\bar4, \bar2) \\
& (1\bar2, \bar3\bar4), (1\bar3, \bar2\bar4), (1\bar4, \bar2\bar3) \\
\end{split}
\end{equation}
and their symmetrical pairs, where $1$ occurs on the right.
The words in these pairs can 
now be decoded into lists of sets partitions according to
the remark following Eq. \eqref{eq:w2lsp}.
This gives the set of Equation \eqref{211}:
the first three words give rise to two pairs of partitions (the
first one gives both $(1\bar2\bar3,\bar4)$ and $(\bar21\bar3,\bar4)$) whereas the
last three ones yield only one.
}
\end{example}

Since we have only used the coproduct rule of the basis
$V_I$, we have as well:

\begin{corollary}\label{cor:generic}
The same is true if one replaces $V_I$ by any basis whose coproduct is the
unshuffle of compositions, for example a basis of product of
primitive elements such as $\Psi^I$ or $\Phi^I$.
\end{corollary}

\subsection{Generating functions and closed formulas}\label{sec:genfunc}

\begin{proposition}
For two compositions $I,J$ of $n$ of the same length $k$,
\begin{equation}
c_{IJ}=\prod_{s=1}^k \binom{j_1+\cdots+j_s-(i_1+\cdots+i_{s-1})-1}{i_s-1}\,.
\end{equation}
\end{proposition}

Examples of the matrices $C=(c_{IJ})$ are given in Section~\ref{sec-mats}.

\medskip
\Proof Let $\pi$ be a set partition such that $C(\pi)=J$, and order
the blocks by increasing maxima. The first block is composed of
$j_1$ and of $i_1-1$ elements strictly smaller than $j_1$, which yields
$\binom{j_1-1}{i_1-1}$ choices. Having chosen those elements, the
second block is composed of $j_1+j_2$, and of $i_2-1$ elements smaller
than $j_1+j_2$, and different from the $i_1$ elements of the first block,
which leaves us with $\binom{j_1+j_2-1-i_1}{i_2-1}$ choices, and so on.
\qed

These expressions have simple generating series, which allow
to find immediately the inverse matrices.

\begin{proposition}
For a composition $J$ of $n$ of length $k$, the generating function
of the column $c_{IJ}$ is
\begin{equation}\label{eqcIJx}
\sum_{I\vDash n,\ \ell(I)=k} c_{IJ}\prod_{s=1}^k x_s^{i_s-1}
=\prod_{s=1}^k (x_s+\cdots +x_k)^{j_s-1}\,. 
\end{equation}
\end{proposition}

\begin{corollary}\label{corZU}
Let $Z_J$ be the dual basis of $X_J$, and $U_I$ be the dual basis
of $V_I$ in $\QSym$. Set
\begin{equation}
Z_I = \sum_J d_{IJ}U_J\,,
\end{equation}
where $D$ is the transpose of $C^{-1}$.
The generating function of row $I$ is
\begin{equation}\label{eqdIJy}
\sum_J d_{IJ}\prod_{s=1}^k y_s^{j_s-1} = \prod_{s=1}^k (y_s-y_{s+1})^{j_s-1}\quad
(\text{with $y_{k+1}:=0$})
\end{equation}
\end{corollary}

\Proof Define $y_i=x_i+\dots+x_k$. Then the system of equations
\begin{equation}
y_1^{j_1-1}y_2^{j_2-1}\cdots y_k^{j_k-1}=
(x_1+\cdots+x_k)^{j_1-1}(x_2+\cdots+x_k)^{j_2-1}\cdots x_k^{j_k-1}
\end{equation} 
for all compositions of $n$ of length $k$ is clearly equivalent to
\begin{equation}\label{eqxy}
x_k=y_k,\quad x_{k-1} = y_{k-1}-y_k,\quad \ldots,\quad x_1=y_1-y_2\,.
\end{equation} 
\qed

For example, we can read on the matrices below that
\begin{equation}
X_{12}=V_{12},\ X_{21}=V_{21}+V_{12},\ \text{so that}\ Z_{21}=U_{21}\
\text{and} \ Z_{12}=U_{12}-U_{21}.
\end{equation}
Hence, since $U_IU_J$ is given by the shuffle of compositions,
\begin{equation}
Z_1Z_{21}= 2U_{211}+U_{121}= Z_{121}+3Z_{211}
\end{equation}
and
\begin{equation}
Z_1Z_{12}= 2U_{211}-2U_{211}= 2Z_{112}+2Z_{121}
\end{equation}
as one can read $U_{211}=Z_{221}$, $U_{121}=Z_{121}+Z_{211}$ and $U_{112}=Z_{112}+Z_{121}+Z_{211}$.
\bigskip

\subsection{Relations with Ecalle's generating functions}

By Corollary \ref{cor:generic},
we can assume that $V_I=Y^I$, where $Y^I=Y_{i_1}\cdots Y_{i_r}$
is the multiplicative basis of $\Sym$ constructed from a generating sequence
$(Y_n)$ of the primitive Lie algebra of $\Sym$. Let $U_I$ be the dual
basis of $Y^I$ in $\QSym$. Then, by Corollary \ref{corZU}, we have
\begin{equation}\label{eqZag}
\sum_{\ell(J)=k}Z_J y_1^{j_1-1}y_2^{j_2-1}\cdots y_k^{j_k-1}
=
\sum_{\ell(I)=k}U_I x_1^{i_1-1}x_2^{i_2-1}\cdots x_k^{i_k-1}\,.
\end{equation}
A similar relation (up to reversal of the indices) occurs in Ecalle's
works on MZVs. If in the l.h.s. we replace $Z_J$ by the integral representation \eqref{eqintZ}
of $\zeta(J)$ (denoted by ${\rm Wa}^J$ in \cite{Ec10}), the l.h.s. becomes what
Ecalle denotes by ${\rm Zag}^{(x_1,\ldots,x_k)}$ (actually, Ecalle works with
colored MZVs and there is a second set of parameters). The equality \eqref{eqZag} translates
into the statement ``${\rm Wa}^\bullet$ symmetral $\Leftrightarrow$ ${\rm Zag}^\bullet$ symmetral''.
Indeed, by definition,  the symmetrality of ${\rm Wa}^\bullet$ is equivalent
to the product formula
\begin{equation}
U_I U_J = \sum_K \< K\ |\ I\shuffle J \> U_K\,,
\end{equation}
which, written in the form of a generating function, reads
\begin{equation}\label{eqZag2}
{\rm Zag}^{(x_1,\ldots,x_k)}{\rm Zag}^{(x_{k+1},\ldots,x_n)}
=\sum_{\sigma\in 12\cdots k\shuffle k+1\cdots n}
{\rm Zag}^{\left(x_{\sigma^{-1}(1)},\ldots,x_{\sigma^{-1}(n)}\right)}.
\end{equation} 


\subsection{Noncommutative Bell polynomials}

To be more specific, one may take the multiplicative basis $Y^I$
constructed from the normalized Dynkin elements $Y_n=(n-1)!\Psi_n$.
Then, 
\begin{equation}
\sum_{J\vDash n}X_J = \sum_{I\vDash n}\beta_I Y^I =  B_n(Y)
\end{equation}
is the noncommutative Bell polynomial in the $Y_i$
as defined in~\cite{SR} by the recursive formula
\begin{equation}
B_{n+1} = \sum_{k=0}^n\binom{n}{k} B_{n-k}Y_{k+1}
\quad\text{and}\quad B_0=1
\end{equation}
which translates into the recursion 
\begin{equation}
(n+1)S_{n+1} = \sum_{k=0}^nS_{n-k}\Psi_{k+1}
\end{equation}
after the change of variables,
so that  for this choice
of the $Y_i$, it reduces to $n!S_n$. 

\subsection{Direct construction from binary words}\label{sec:binw}

Rather surprisingly, the coefficients $c_{IJ}$, which have been obtained
in a canonical way from $\WSym$ without choosing a basis
of $\Sym$, can also be obtained by a very specific choice, the above
normalized version of the $\Psi$ basis, interpreted in terms
of Lie polynomials in two letters. 

Dual to the realization of $\QSym$ as a subalgebra of $\KK_{\shuffle}\<a,b\>$,
there is a simple realization of $\Sym$ as a subalgebra of $\KK\<a,b\>$,
regarded as the universal enveloping algebra of $L(a,b)$, the free
Lie algebra on two letters. If one sets
\begin{equation}
\Psi_n = \frac1{(n-1)!} \ad_a^{n-1}b = \frac1{(n-1)!} [a,[a,[\ldots,[a,b]\ldots]]]
\end{equation}
then there is a simple expression for $S_n$:
\begin{equation}
S_n = \frac1{n!}\sum_{k=0}^n\binom{n}{k}(-1)^{k}(a+b)^{n-k}a^k\,.
\end{equation}
Indeed, 
\begin{equation}
\psi(t):=\sum_{n\ge 1}t^{n-1}\Psi_n = e^{t\ad_a}b = e^{ta}b e^{-ta}
\end{equation}
so that the differential equation for $x(t)=\sum_n t^n S_n$, which is
\begin{equation}
x'(t)=x(t)\psi(t)\quad \text{with $x(0)=1$}
\end{equation}
yields
\begin{equation}
x(t)=e^{t(a+b)}e^{-ta}\,.
\end{equation}

The algebra morphism defined by
$\iota: (n-1)!\Psi_n\mapsto Y_n=\ad^{n-1}_ab$ is also
a coalgebra morphism, as it sends primitive elements to primitive elements.
However, since the $Z$-bases define embeddings of $\QSym$ as a subalgebra
of $\KK_{\shuffle}\<a,b\>b$, by duality, $\Sym$ should be a quotient
coalgebra of $\KK\<a,b\>b$. 

Let $\p:\ \KK\<a,b\>\rightarrow \KK\<a,b\>b$ be the projection defined
by $\p(w)=0$ if $w$ ends by $a$, and let 
\begin{equation}
\Delta' = (\p\otimes\p)\circ \Delta
\end{equation}
where $\Delta$ is the usual coproduct of  $\KK\<a,b\>$ (for which the letters
are primitive). Then, 
\begin{equation}
\p(Y_n)=a^{n-1}b = W_n 
\end{equation}
is primitive for $\Delta'$. For $I=(i_1,\ldots,i_r)$,
\begin{equation}
\p (Y^I) = Y^{i_1,\ldots,i_{r-1}}\p(Y_{i_r})
\end{equation}
and by induction on $r$,
\begin{equation}
\Delta'\circ\p (Y^I) = (\p\otimes \p)\circ \Delta (Y^I)
\end{equation}
so that $\p$ induces an isomorphism of coalgebras
\begin{equation}
\p:\ (\iota(\Sym),\Delta)\longrightarrow (\KK\<a,b\>b,\Delta')\,.
\end{equation}
Now,
\begin{equation}
Y_n = \sum_{k=0}^{n-1}(-1)^k\binom{n-1}{k} a^{n-1-k}ba^k
\end{equation}
and a straightforward calculation 
(see \cite[Appendix A]{Rac}) shows that
\begin{equation}
\p(Y^I) =\sum_{J}d_{IJ}W_J
\end{equation}
where the $d_{IJ}$ are as in Corollary \ref{corZU}.
By duality, the image $Z_I$ in $\QSym$ of the dual basis $W_I^*$
by $\iota^*\circ\p^*$ satisfies the binary shuffle
relations \eqref{shufZ}. 
These considerations, which are essentially in Racinet's thesis \cite{Rac}
explain the appearance of the generating functions \eqref{eqdIJy} and
\eqref{eqcIJx} in the theory of MZVs, and in particular Ecalle's swap
operation. We can now see that they follow from a particular choice of a
generating sequence of the primitive Lie algebra of $\Sym$, together with a
specific embedding of $\Sym$ in $\KK\<a,b\>$.

\medskip
Finally, it is instructive to play a little with the noncommutative
Bell polynomials in this context.
Setting $Y_k=\ad_a^{k-1}b$ in the noncommutative Bell polynomials,
this yields
\begin{equation}\label{Snab}
B_n = \sum_{k=0}^n\binom{n}{k}(-1)^{k}(a+b)^{n-k}a^k\,.
\end{equation}
Let now $L$ be the linear operator
\begin{equation}
L =\ad_a+b\,.
\end{equation}
Let $D$ denote the derivation $\ad_a$,
and set $b^{(k)}=D^k b$.
According to \cite[Theorem 2]{SR},
\begin{equation}
B_n(b,b',b'',\ldots,b^{(n-1)})  = L^n (1).
\end{equation}
Thus, the
above considerations amount to that fact
$L^n(1)$ is equal to the r.h.s of \eqref{Snab},
which can of course be proved directly:
\begin{equation}
x(t) := e^{t(a+b)}e^{-ta}
\end{equation}
satisfies obviously
\begin{equation}
\frac{dx}{dt} = (a+b)x - xa = Lx,\quad x(0)=1\,.
\end{equation}

Comparing~\eqref{eqcIJx} with~\cite[Theorem 2]{SR}, we obtain the
binomial identity
\begin{corollary}
The column sums of the matrices $(c_{IJ})$ are
\begin{equation}
\sum_{J\vDash n,\ \ell(J)=k}
  \prod_{s=1}^k \binom{j_1+\cdots+j_s-(i_1+\cdots+i_{s-1})-1}{i_s-1}
=
\prod_{s=2}^k\binom{i_1+\cdots+i_s-1}{i_1+\cdots+i_{s-1}}\,.
\end{equation}
\end{corollary}

\section{Matrices}
\label{sec-mats}

The entry in row $I$ and column $J$ is the coefficient of $V_I$ in $X_J$ in
the first matrix, and the coefficient of $U_I$ in $Z_J$ in the second one.

{Case $n=3$, $k=2$}
\begin{equation}
\left(\begin{array}{rr}
1 & . \\
1 & 1
\end{array}\right)
\qquad
\left(\begin{array}{rr}
1 & -1 \\
. & 1
\end{array}\right)
\end{equation}

Case $n=4$ with $k=2$ and $k=3$.
\begin{equation}
\left(\begin{array}{rrr}
1 & . & . \\
2 & 1 & . \\
1 & 1 & 1
\end{array}\right)
\qquad
\qquad
\left(\begin{array}{rrr}
1 & -2 & 1 \\
. & 1 & -1 \\
. & . & 1
\end{array}\right)
\end{equation}

\begin{equation}
\left(\begin{array}{rrr}
1 & . & . \\
1 & 1 & . \\
1 & 1 & 1
\end{array}\right)
\qquad
\qquad
\left(\begin{array}{rrr}
1 & -1 & . \\
. & 1 & -1 \\
. & . & 1
\end{array}\right)
\end{equation}

Case $n=5$, for all values of $k$ from $2$ to $4$.
\begin{equation}
\left(\begin{array}{rrrr}
1 & . & . & . \\
3 & 1 & . & . \\
3 & 2 & 1 & . \\
1 & 1 & 1 & 1
\end{array}\right)
\qquad
\qquad
\left(\begin{array}{rrrr}
1 & -3 & 3 & -1 \\
. & 1 & -2 & 1 \\
. & . & 1 & -1 \\
. & . & . & 1
\end{array}\right)
\end{equation}

\begin{equation}
\left(\begin{array}{rrrrrr}
1 & . & . & . & . & . \\
2 & 1 & . & . & . & . \\
2 & 1 & 1 & . & . & . \\
1 & 1 & . & 1 & . & . \\
2 & 2 & 1 & 2 & 1 & . \\
1 & 1 & 1 & 1 & 1 & 1
\end{array}\right)
\quad
\left(\begin{array}{rrrrrr}
1 & -2 & . & 1 & . & . \\
. & 1 & -1 & -1 & 1 & . \\
. & . & 1 & . & -1 & . \\
. & . & . & 1 & -2 & 1 \\
. & . & . & . & 1 & -1 \\
. & . & . & . & . & 1
\end{array}\right)
\end{equation}

\begin{equation}
\left(\begin{array}{rrrr}
1 & . & . & . \\
1 & 1 & . & . \\
1 & 1 & 1 & . \\
1 & 1 & 1 & 1
\end{array}\right)
\qquad
\qquad
\left(\begin{array}{rrrr}
1 & -1 & . & . \\
. & 1 & -1 & . \\
. & . & 1 & -1 \\
. & . & . & 1
\end{array}\right)
\end{equation}

For $k=2$, the block is always given by the Pascal triangle.
For $k=n-1$, it is the lower triangular matrix with all entries equal
to~1. Here follows the remaining blocks for $n=6$, that is, for $k=3$ and $4$.

{\small
\begin{equation*}
\left(\begin{array}{rrrrrrrrrr}
1 & . & . & . & . & . & . & . & . & . \\
3 & 1 & . & . & . & . & . & . & . & . \\
3 & 1 & 1 & . & . & . & . & . & . & . \\
3 & 2 & . & 1 & . & . & . & . & . & . \\
6 & 4 & 2 & 2 & 1 & . & . & . & . & . \\
3 & 2 & 2 & 1 & 1 & 1 & . & . & . & . \\
1 & 1 & . & 1 & . & . & 1 & . & . & . \\
3 & 3 & 1 & 3 & 1 & . & 3 & 1 & . & . \\
3 & 3 & 2 & 3 & 2 & 1 & 3 & 2 & 1 & . \\
1 & 1 & 1 & 1 & 1 & 1 & 1 & 1 & 1 & 1
\end{array}\right)
\quad
\left(\begin{array}{rrrrrrrrrr}
1 & -3 & . & 3 & . & . & -1 & . & . & . \\
. & 1 & -1 & -2 & 2 & . & 1 & -1 & . & . \\
. & . & 1 & . & -2 & . & . & 1 & . & . \\
. & . & . & 1 & -2 & 1 & -1 & 2 & -1 & . \\
. & . & . & . & 1 & -1 & . & -1 & 1 & . \\
. & . & . & . & . & 1 & . & . & -1 & . \\
. & . & . & . & . & . & 1 & -3 & 3 & -1 \\
. & . & . & . & . & . & . & 1 & -2 & 1 \\
. & . & . & . & . & . & . & . & 1 & -1 \\
. & . & . & . & . & . & . & . & . & 1
\end{array}\right)
\end{equation*}
}

{\small
\begin{equation*}
\left(\begin{array}{rrrrrrrrrr}
1 & . & . & . & . & . & . & . & . & . \\
2 & 1 & . & . & . & . & . & . & . & . \\
2 & 1 & 1 & . & . & . & . & . & . & . \\
2 & 1 & 1 & 1 & . & . & . & . & . & . \\
1 & 1 & . & . & 1 & . & . & . & . & . \\
2 & 2 & 1 & . & 2 & 1 & . & . & . & . \\
2 & 2 & 1 & 1 & 2 & 1 & 1 & . & . & . \\
1 & 1 & 1 & . & 1 & 1 & . & 1 & . & . \\
2 & 2 & 2 & 1 & 2 & 2 & 1 & 2 & 1 & . \\
1 & 1 & 1 & 1 & 1 & 1 & 1 & 1 & 1 & 1
\end{array}\right)
\quad
\left(\begin{array}{rrrrrrrrrr}
1 & -2 & . & . & 1 & . & . & . & . & . \\
. & 1 & -1 & . & -1 & 1 & . & . & . & . \\
. & . & 1 & -1 & . & -1 & 1 & . & . & . \\
. & . & . & 1 & . & . & -1 & . & . & . \\
. & . & . & . & 1 & -2 & . & 1 & . & . \\
. & . & . & . & . & 1 & -1 & -1 & 1 & . \\
. & . & . & . & . & . & 1 & . & -1 & . \\
. & . & . & . & . & . & . & 1 & -2 & 1 \\
. & . & . & . & . & . & . & . & 1 & -1 \\
. & . & . & . & . & . & . & . & . & 1
\end{array}\right)
\end{equation*}
}

Comparing these matrices with the matrices given in Equations~\eqref{n4k23}
and~\eqref{n5k23}, on sees that the transposes of the matrices giving
the coefficient of $Z_I$ on $U_J$ has many similarities with those 
constructed by the recursive procedure of Section\ref{sec:algo}.

\section{Rational functions and colored free quasi-symmetric functions}

Recall that $\QSym$ is the commutative image of the algebra $\FQSym$
of free quasi-symmetric functions, that is, $F_I(X)$ is obtained
by sending the noncommutating variables $a_i$ of
\begin{equation}
\F_\sigma(A)=\sum_{\std(w)=\sigma^{-1}}w
\end{equation}
to commuting variables $x_i$. 

In \cite{CHNT}, the vector space $\FQSym$ is identified to the Zinbiel
operad, realized as a suboperad of the operad of rational moulds.

A {\em mould}, as defined by Ecalle, is a ``function of a variable
number of variables'', that is, a sequence $f=(f_n(u_1,\ldots,u_n))$ of
functions of $n$ (continuous or discrete) variables. There is a bilinear
operation defined componentwise by
\begin{equation}
{\rm mu}(f_n,g_m) =f_n(u_1,\ldots,u_n)g_m(u_{n+1},\ldots,u_{n+m})
\end{equation}
which defines an associative product $*$ on homogeneous moulds (those with
only one $f_n$ nonzero). For this product, the rational functions
\begin{equation}
f_\sigma(u_1,\ldots,u_n)=\frac1{u_{\sigma(1)}(u_{\sigma(1)}+u_{\sigma(2)})\cdots
(u_{\sigma(1)}+u_{\sigma(2)}+\cdots+u_{\sigma(n)})}
\end{equation}
span a subalgebra isomorphic to $\FQSym$ under the correspondence
$f_\sigma\mapsto \F_\sigma$.

In \cite{GX}, it is proved that the fractions
\begin{equation}
z_{\sigma,s}(u_1,\ldots,u_n)=\frac1{u_{\sigma(1)}^{s_1}(u_{\sigma(1)}+u_{\sigma(2)})^{s_2}\cdots
(u_{\sigma(1)}+u_{\sigma(2)}+\cdots+u_{\sigma(n)})^{s_n}}
\end{equation}
where $\sigma\in\SG_n$ and $s\in\NN^{*n}$ satisfy a product formula generalizing the
binary shuffle \eqref{shufZ}. That is, if one sets
\begin{equation}
\epsilon_{\sigma,s} = 0^{s_1-1}\sigma_1 0^{s_2-1}\sigma_2 \cdots 0^{s_n-1}\sigma_n
\end{equation}
then,
\begin{equation}
z_{\sigma,s}*z_{\tau,t}
 = \sum_{\omega,w}
   \<\epsilon_{\omega,w}| \epsilon_{\sigma,s}\shuffle \epsilon_{\tau[n],t}\>
      z_{\omega,w}
\end{equation}
where $\tau[n]$ denotes as usual the permutation $\tau$ shifted by the length
of $\sigma\in\SG_n$.
There is a natural bigrading $\deg z_{\sigma,s} = (n,|s|)$ if $\sigma\in\SG_n$
and $|s|=\sum s_i$.

Again, there exist colored versions of $\FQSym$, based on symbols
$\F_{\sigma,c}$ where $\sigma\in\SG_n$ and $c\in C^n$ is a color word, $C$
being the color alphabet.
Taking $C=\NN^*$, we obtain a bigraded vector space isomorphic to the linear
span of the $z_{\sigma,s}$.
The product rule is
\begin{equation}
\F_{\sigma',c'}\F_{\sigma'',c''}
=
\sum_{(\sigma,c)\in (\sigma',c')\Cup (\sigma'',c'')}\F_{\sigma,c}\,,
\end{equation}
where the shifted shuffle $\Cup$ of colored permutations is 
computed by shifting the letters of $\sigma''$ by the size of
$\sigma'$ and shuffling the colored letters.

Thus, $\F_{\sigma,s}\mapsto z_{\sigma,s}$ is not an
isomorphism of algebras, but 
the previous considerations allow us to prove:

\begin{theorem} 
Define a basis $\Z_{\sigma,J}$ by
\begin{equation}
\F_{\sigma, I}= \sum_{|J|=|I|}c_{IJ}\Z_{\sigma,J}
\end{equation}
for $\sigma\in\SG_n$ and $I$ a composition of length $n$, the $c_{IJ}$
being as in \eqref{eqcIJx}. Then, 
\begin{equation}
\Z_{\sigma',J'}\Z_{\sigma'',J''}
 = 
\sum_{\sigma,J}\<\epsilon_{\sigma,J}| \epsilon_{\sigma',I'}\shuffle \epsilon_{\sigma''[n],I''}\>
\Z_{\sigma,J}.
\end{equation}
\end{theorem}

\Proof
This can be seen, for example, by a straightforward generalization
of the argument of Section \ref{sec:binw}. Replace the letter $b$
by an alphabet $B=\{b_i|i\ge 1\}$, and define elements of $\KK\<a,B\>$
by
\begin{equation}
Y^{i,n} = \ad_a^{n-1}b_i,\quad Y^{u,I}=Y^{u_1,i_1}\cdots Y^{u_r,i_r}
\end{equation}
for a  word $u=u_1\cdots u_r$ and a composition $I=(i_1,\ldots,i_r)$.

Let, as before $\p:\ \KK\<a,B\>\rightarrow \bigoplus_ {b\in B}\KK\<a,B\>b$ 
be the projection defined
by $\p(w)=0$ if $w$ ends by $a$, and  
\begin{equation}
\Delta' = (\p\otimes\p)\circ \Delta
\end{equation}
where $\Delta$ is the usual coproduct of  $\KK\<a,B\>$ for which the letters
are primitive. Then, 
\begin{equation}
\p(Y^{i,n})=a^{n-1}b_i =: W_{i,n} 
\end{equation}
is primitive for $\Delta'$,
and
\begin{equation}
\p (Y^{u,I}) = Y^{u_1,i_1}\cdots Y^{u_{r-1},i_{r-1}}\p(Y^{u_r,i_r})
\end{equation}
so that 
\begin{equation}
\Delta'\circ\p (Y^{u,I}) = (\p\otimes \p)\circ \Delta (Y^{u,I})
\end{equation}
Assuming now that $u=\sigma$ is a permutation, 
and taking into account the product rule of
$\FQSym^{(\NN^*)}$,
we can also write
this coproduct in the form
\begin{equation}\label{eq:deltacol}
\Delta'\circ\p (Y^{u,I}) = 
\sum_{u',I';\,u'',I''}
\< \F_{\sigma,I} | \F_{\std(u'),I'}\F_{\std(u''),I''} \>
Y^{u',I'}\otimes Y^{u'',I''}.
\end{equation}
Again, 
\begin{equation}
\p(Y^{\sigma,I}) =\sum_{J}d_{IJ}W_{\sigma,J},
\end{equation}
and the result follows from \eqref{eq:deltacol} by duality. 
\qed

Taking into account the generating functions of Section \ref{sec:genfunc}, we
have
\begin{equation}
\sum_{\ell(I)=n}y_1^{i_1-1}y_2^{i_2-1}\cdots y_n^{i_n-1}
\Z_{\sigma,I}
=
\sum_{\ell(I)=n}
(y_1-y_2)^{i_1-1}(y_2-y_3)^{i_2-1}\cdots y_n^{i_n-1}
\F_{\sigma,I}\,.
\end{equation}

Finally, it can be shown that the vector space spanned by the $z_{\sigma,I}$
is stable under the operadic compositions of \cite{CHNT}. The resulting
operad will be investigated in a separate paper.

\section*{Acknowledgements}
This research has been partially supported by the project CARMA (ANR-12-BS01-0017)
of the Agence Nationale de la Recherche and the Labex B\'ezout of Universit\'e Paris-Est.


\end{document}